\documentclass[12pt,leqno]{article}
\usepackage{latexsym}
\usepackage{bbm}
\usepackage{amssymb,amsmath}
\begin{titlepage}
\title{The large sieve with sparse sets of moduli}
\author{Stephan Baier}
\date{22.08.05}
\end{titlepage}
\begin{document}
\maketitle
$ $\\
$ $\\
{\bf Address of the author:}\medskip\\ 
Stephan Baier\\ 
Jeffery Hall\\ 
Department of Mathematics and Statistics\\ 
Queen`s University\\
University Ave\\
Kingston, Ontario, Canada\\
K7L 3N6\medskip\\ 
e-mail: sbaier@mast.queensu.ca
\newpage 
$ $\\
{\bf Abstract:} Extending a method of D. Wolke \cite{Wol}, 
we establish a general 
result on the large sieve with sparse
sets ${\cal{S}}$ of moduli which are in a sense well-distributed in 
arithmetic 
progressions. We then apply our result to the case when ${\cal{S}}$ consists
of sqares. In this case we obtain an estimate which improves a recent result by L. Zhao \cite{Zha}.\\ \\
Mathematics Subject Classification (2000): 11N35, 11L07, 11B57\\  \\
Key words: large sieve, sparse set of moduli, square moduli, Farey fractions
in short intervals, estimates on exponential sums\\ \\ 
\section{Introduction}
Throughout this paper, we reserve the symbols 
$c$, $c_i$ $(i=1,2,...)$ for absolute 
constants. Further, we suppose that 
$(a_n)$ is a sequence of complex numbers and that $Q,N\ge 1$. We set
\begin{equation}
S(\alpha):=\sum\limits_{n\le N} a_n e(n\alpha).\label{0}
\end{equation}

Bombieri's \cite{Cla} classical large sieve inequality asserts that 
\begin{equation}
\sum\limits_{q\le Q} \sum\limits_{\scriptsize \begin{array}{cccc} 
a=1\\(a,q)=1\end{array}}^q \left\vert S\left(\frac{a}{q}\right)\right\vert^2
\ll (N+Q^2)Z,\label{1}
\end{equation}
where 
$$
Z:=\sum\limits_{n\le N} \vert a_n\vert^2.
$$
It is natural to ask whether (\ref{1}) can be improved if 
the moduli $q$ run over a sparse set ${\cal{S}}$ of natural numbers $\le Q$.
In the sequel, let $S$ be the cardinality of ${\cal{S}}$.  

P.D.T.A. Elliott \cite{Ell} considered the case when ${\cal{S}}$ consists of 
primes. 
Essentially, he proved that 
\begin{equation}
\sum\limits_{p\in {\cal{S}}} \sum\limits_{a=1}^{p-1} 
\left\vert S\left(\frac{a}{p}\right)\right\vert^2
\ll \left(\frac{N^2}{\log N}+QS\right)Z\label{A}
\end{equation}
and conjectured that the right-hand side can always 
be replaced by 
\begin{equation}
\ll (N+QS)Z. \label{B}
\end{equation}
 
For the case when ${\cal{S}}$ consists of all primes
$p\le Q$ D. Wolke \cite{Wol} found the estimate 
\begin{equation}
\sum\limits_{p\le Q} \sum\limits_{a=1}^{p-1} 
\left\vert S\left(\frac{a}{p}\right)\right\vert^2
\le \frac{c}{1-\delta} \frac{Q^2\log\log Q}{\log Q}Z \label{2}
\end{equation}
if $Q\ge 10$, $N=Q^{1+\delta}$, $0<\delta<1$. 

L. Zhao \cite{Zha} studied the case when ${\cal{S}}$ consists of squares. 
He established the bound
\begin{equation}
\sum\limits_{q\le Q_1} \sum\limits_{\scriptsize \begin{array}{cccc} 
a=1\\(a,q)=1\end{array}}^{q^2}
\left\vert S\left(\frac{a}{q^2}\right)\right\vert^2
\ll (\log 2Q_1)\left(Q_1^3+(N\sqrt{Q_1}+\sqrt{N}Q_1^2)N^{\varepsilon}\right)Z. 
\label{3}
\end{equation}
The classical form (\ref{1}) of the large sieve implies only the bound
\begin{equation}
\sum\limits_{q\le Q_1} \sum\limits_{\scriptsize \begin{array}{cccc} 
a=1\\(a,q)=1\end{array}}^{q^2}
\left\vert S\left(\frac{a}{q^2}\right)\right\vert^2
\ll (N+Q_1^{4})Z, \label{Z1}
\end{equation}
which is weaker than (\ref{3}) if $Q_1\gg N^{2/7+\varepsilon}$. 
Using the later Lemma 1,
which gives a general large sieve bound, 
one can also show that (see \cite{Zha})
\begin{equation}
\sum\limits_{q\le Q_1} \sum\limits_{\scriptsize \begin{array}{cccc} 
a=1\\(a,q)=1\end{array}}^{q^2}
\left\vert S\left(\frac{a}{q^2}\right)\right\vert^2
\ll Q_1(N+Q_1^{2})Z, \label{Z2}
\end{equation}
which is weaker than (\ref{3}) if $Q_1\ll N^{1/2-\varepsilon}$. 
So (\ref{3}) is sharper
than both (\ref{Z1}) and (\ref{Z2}) if $N^{2/7+\varepsilon}\ll Q_1\ll 
N^{1/2-\varepsilon}$. 

The aim of the present paper is to establish a general large sieve bound 
for the case when ${\cal{S}}$ is a sparse set which is in a
certain sense well-distributed in arithmetic progressions. 
To do so, we shall use an extension of   
the method developed in \cite{Wol}. Our main result, the later Theorem 2, will
imply Wolke's estimate $(\ref{2})$ and a sharpened version of Zhao's estimate
$(\ref{3})$. For the sake of its generality, the later Theorem 2
might have a number of
other applications.\\
  
\section{Statement of the results}
In the sequel, for $t\in \mathbbm{N}$ let
$$
{\cal{S}}_t:=\{q\in \mathbbm{N}\ :\ tq\in {\cal{S}}\}
$$
and $S_t:=\vert {\cal{S}}_t \vert$.   
We first give a general estimate for the sum in question in terms of the number
of elements of ${\cal{S}}_t$ in short segments of arithmetic progressions.\\
 
{\bf Theorem 1:} \begin{it} Suppose that $0\le M\le Q$ and  
${\cal{S}}\subset (M,M+Q]$. 
For $u>0$, $t,k\in \mathbbm{N}$, $l\in\mathbbm{Z}$ define 
$$
A_t(u,k,l):=\max\limits_{M/t\le y\le (M+Q)/t} 
\vert \{q\in {\cal{S}}_t\cap (y,y+u] \ : \ q\equiv l\mbox{ mod }k\} \vert.
$$
Put 
$$
U:=\left\{\begin{array}{llll} 1, & \mbox{ if }\ M< \sqrt{N},\\ \\
0, & \mbox{ otherwise.}\end{array}\right.
$$ 
Then,
\begin{eqnarray}
\label{4} & &\\
\sum\limits_{q\in {\cal{S}}} \sum\limits_{\scriptsize \begin{array}{cccc} 
a=1\\(a,q)=1\end{array}}^q \left\vert S\left(\frac{a}{q}\right)\right\vert^2
&\le& c_1NZ\left(U+\max\limits_{r\le \sqrt{N}}\ 
\max\limits_{1/N\le z\le 1/(r\sqrt{N})} \max\limits_{\scriptsize
\begin{array}{cccc} h\in \mathbbm{Z}\\(h,r)=1\end{array}}\right. \nonumber\\
& & \left. \sum\limits_{t\vert r}
\sum\limits_{\scriptsize \begin{array}{cccc} 
0<m\le 4rzQ/t\\ (m,r/t)=1\end{array}} 
A_t\left(\frac{Q}{tzN},\frac{r}{t},hm\right)\right).\nonumber
\end{eqnarray} 
\end{it} 

To simplify the estimate (\ref{4}), we now set some natural conditions to the 
terms $A_t(u,k,l)$ with $(k,l)=1$. We motivate these conditions by 
heuristic ideas.

Suppose that $t\le \sqrt{N}$ and $0<u\le Q/t$.
If the set ${\cal{S}}_t$ is nearly evenly distributed on the interval 
$(M/t,(M+Q)/t]$, we would expect that
$$
A_t(u,1,0)\le C\left(1+\frac{S_t}{Q/t}\cdot u\right)
$$
with $C\ge 1$ not too large.   
If we, more generally, 
assume the set ${\cal{S}}_t$ to be nearly evenly distributed in the residue
classes mod $k$, we get 
\begin{equation}
A_t(u,k,l)\le C\left(1+\frac{S_t/k}{Q/t}\cdot u\right). \label{20}
\end{equation}
For many sets ${\cal{S}}_t$ naturally appearing in arithmetic, 
the term on the right side of (\ref{20}) 
will roughly give the 
correct order of magnitude, but in addition we have to take 
into consideration possible fluctuations on varying $l$. 
Therefore we multiply the right side of (\ref{20}) by a non-negative constant
$\delta_t(k,l)$ depending on the residue class $l$ mod $k$.
Thus, our first condition to $A_t(u,k,l)$ is
\begin{equation}
A_t(u,k,l)\le C\left(1+\frac{S_t/k}{Q/t}\cdot u\right)\delta_t(k,l).\label{5}
\end{equation}

Since $\sum_{l=1}^{k} A_t(u,k,l)$
should be roughly of the same size as $A_t(u,1,0)$, we may further assume that
\begin{equation}
\sum\limits_{\scriptsize\begin{array}{cccc} l=1\\(k,l)=1\end{array}}^{k} 
\delta_t(k,l)\le k, \label{6}
\end{equation}
which is our second condition. 

Finally, we suppose that 
\begin{equation}
\delta_t(k,l)\le X \label{7}
\end{equation}
for all $k\le \sqrt{N}/t$, $l\in \mathbbm{Z}$ with $(k,l)=1$, 
where we think of $X$ as a 
quantity which is small compared with $Q$ and $N$. 

In the sections 3 and 4 we will examplify our heuristics by the cases when
${\cal{S}}$ consists of primes as well as of squares. 

We shall later bound the right side of (\ref{4}) by using the conditions 
(\ref{5}), (\ref{6}), (\ref{7}). This shall lead us to the following\\

{\bf Theorem 2:} \begin{it}  
Suppose that $0\le M\le Q$ and ${\cal{S}}\subset (M,M+Q]$. 
Define $A_t(u,k,l)$ and $U$ as in Theorem 1. 
Assume the conditions (\ref{5}), (\ref{6}), (\ref{7}) to hold for all   
$t\le \sqrt{N}$, $k\le \sqrt{N}/t$, $l\in \mathbbm{Z}$ with $(k,l)=1$, and  
$kQ/\sqrt{N}\le u\le Q/t$. Then,
\begin{eqnarray}
\label{8} & &\\
& & \sum\limits_{q\in {\cal{S}}} \sum\limits_{\scriptsize \begin{array}{cccc} 
a=1\\(a,q)=1\end{array}}^q \left\vert S\left(\frac{a}{q}\right)\right\vert^2
\nonumber\\ &\le& c_2C\left(NU+ 
(\min\{QX,N\}+Q)\left(\sqrt{N}\log\log 10N+ 
\max\limits_{r\le\sqrt{N}} \sum\limits_{t\vert r} S_t\right)\right)Z. 
\nonumber
\end{eqnarray} 
\end{it} 

Using $S_t\le S$ for all $t\in \mathbbm{N}$ and $\sum\limits_{t\vert r} 1 \ll
r^{\varepsilon}$, we deduce from Theorem 2\\

{\bf Corollary:} \begin{it}  
Let the conditions and assumptions of Theorem 2 be kept. Fix any 
$\varepsilon>0$. Then,
$$ 
\sum\limits_{q\in {\cal{S}}} \sum\limits_{\scriptsize \begin{array}{cccc} 
a=1\\(a,q)=1\end{array}}^q \left\vert S\left(\frac{a}{q}\right)\right\vert^2
\ll \left(N+QN^{\varepsilon}X\left(\sqrt{N}+S\right)\right)Z.
$$
\end{it} 

This bound should be compared with Elliot's estimate (\ref{A}) and conjecture
(\ref{B}) for the case when 
${\cal{S}}$ consists of primes.
   
Employing Theorem 2 with ${\cal{S}}$ the set of squares lying in the
interval $(Q,2Q]$, we shall 
later derive the following improvement of Zhao's bound
(\ref{3}) from Theorem 2.\\

{\bf Theorem 3:} \begin{it} Suppose that  $Q_1\ge 1$. Fix any 
$\varepsilon>0$. Then,
\begin{equation}
\sum\limits_{q\le Q_1} \sum\limits_{\scriptsize \begin{array}{cccc} 
a=1\\(a,q)=1\end{array}}^{q^2}
\left\vert S\left(\frac{a}{q^2}\right)\right\vert^2
\ll (\log\log 10NQ)^2(Q_1^{3}+N+N^{1/2+\varepsilon}Q_1^2)Z. \label{9}
\end{equation}\end{it}

The bound (\ref{9}) is sharper than the three bounds (\ref{3}), (\ref{Z1})
and (\ref{Z2}) if 
$N^{1/4+\varepsilon}\ll Q_1\ll N^{1/3-\varepsilon}$. \\

\section{Theorem 2 implies Wolke's bound (\ref{2})}
In this section we shall see that Theorem 2 contains Wolke's bound (\ref{2})
as a special case. 

As in \cite{Wol}, we assume that $Q\ge 10$, 
$N=Q^{1+\delta}$, $0<\delta<1$. Further, we choose
$M:=0$ and ${\cal{S}}$ to be the set of all primes $\le Q$. 

First we consider the case when $t=1$. Then, by the 
Brun-Titchmarsh inequality (see \cite{MVa}), we have 
\begin{equation}
A_1(u,k,l)\le \frac{2u}{k\log u/k}\cdot \frac{k}{\varphi(k)} \label{50}
\end{equation}
if $u\ge k$.

For $u\ge kQ/\sqrt{N}$ the inequality (\ref{50}), our assumption 
$N=Q^{1+\delta}$ and the prime number theorem yield 
$$
A_1(u,k,l)\le \frac{4u}{(1-\delta)k\log Q}\cdot \frac{k}{\varphi(k)}
\le \frac{c_3uS}{(1-\delta)kQ} \cdot \frac{k}{\varphi(k)}.
$$
Thus, the condition (\ref{5}) is satisfied
if we set $C:=c_3/(1-\delta)$ and $\delta_1(k,l):=k/\varphi(k)$.
Obviously, the condition (\ref{6}) is also satisfied, and the condition 
(\ref{7}) holds with $X=c_4\log\log Q$ (see \cite{Han}). 

Now, we assume that $t\ge 2$. Then, $S_t\le 1$. Thus, the conditions 
(\ref{5}), (\ref{6}), (\ref{7}) hold trivially with $\delta_t(k,l)=1$ and
$C,X$ chosen as above. 

Therefore, we can apply Theorem 2 to this situation. Taking into account
$S\le c_5Q/\log Q$ by the prime number theorem, we obtain the 
bound (\ref{2}) from (\ref{8}). $\Box$\\

\section{Proof of Theorem 3}
Next, we shall derive Theorem 3 from Theorem 2.
 
If $Q_1^4\le N$, then 
(\ref{9}) follows from (\ref{Z1}). In the following, we assume that
$Q_1^4>N$.

First, we rewrite the sum in question in the form
\begin{equation}
\sum\limits_{q\le Q_1} \sum\limits_{\scriptsize \begin{array}{cccc} 
a=1\\(a,q)=1\end{array}}^{q^2}
\left\vert S\left(\frac{a}{q^2}\right)\right\vert^2
= \sum\limits_{q\in {\cal{S}}} \sum\limits_{\scriptsize \begin{array}{cccc} 
a=1\\(a,q)=1\end{array}}^q
\left\vert S\left(\frac{a}{q}\right)\right\vert^2, \label{21}
\end{equation}
where ${\cal{S}}$ is the set of squares $\le Q_1^2$. 
We spit up the set ${\cal{S}}\cap (\sqrt{N},Q_1^2]$ into 
subsets of the form
$$
{\cal{S}}(Q):={\cal{S}}\cap (Q,2Q],
$$
where $Q\ge \sqrt{N}$.
Our aim is to estimate the partial sums
$$
\sum\limits_{q\in {\cal{S}}(Q)} \sum\limits_{\scriptsize \begin{array}{cccc} 
a=1\\(a,q)=1\end{array}}^q
\left\vert S\left(\frac{a}{q}\right)\right\vert^2.
$$
For the remaining sum we have 
\begin{equation}
\sum\limits_{q\in {\cal{S}}\cap (0,\sqrt{N}]} 
\sum\limits_{\scriptsize \begin{array}{cccc} 
a=1\\(a,q)=1\end{array}}^q
\left\vert S\left(\frac{a}{q}\right)\right\vert^2 \ll NZ \label{22}
\end{equation}
from (\ref{1}).

As previously, we define 
$$
{\cal{S}}_t(Q):=\{q\in \mathbbm{N}\ :\ tq\in {\cal{S}}(Q)\}
$$
and $S_t(Q):=\vert {\cal{S}}_t(Q) \vert$. 
We now determine the set ${\cal{S}}_t(Q)$. Let $t=p_1^{v_1}\cdots p_n^{v_n}$
be the prime number factorization of $t$. For $i=1,...,n$ let
$$
u_i:=\left\{\begin{array}{llll} v_i, & \mbox{ if } v_i \mbox{ is even,}\\ \\
v_i+1, & \mbox{ if } v_i \mbox{ is odd.} \end{array}\right.
$$
Put 
$$
f_t:=p_1^{u_1/2}\cdots p_n^{u_n/2}.
$$
Then $q=q_1^2\in{\cal{S}}$ is divisible by $t$ iff $q_1$ is divisible by
$f_t$. Thus, 
$$
{\cal{S}}_t(Q)=\left\{q_2^2g_t\ : \ \sqrt{Q}/f_t<q_2\le \sqrt{2Q}/f_t\right\}
\subset (Q/t,2Q/t],
$$
where  
$$
g_t:=\frac{f_t^2}{t}=p_1^{u_1-v_1}\cdots p_n^{u_n-v_n}.
$$
Hence,
\begin{equation}
{S}_t(Q)\le \frac{\sqrt{2Q}}{f_t}.\label{25}
\end{equation}

As previously, we suppose that 
$k\in \mathbbm{N}$, $l\in \mathbbm{Z}$ and $(k,l)=1$, and 
define 
$$
A_t(u,k,l):=\max\limits_{Q/t\le y\le 2Q/t} 
\vert \{q\in {\cal{S}}_t(Q)\cap (y,y+u] \ : \ q\equiv l\mbox{ mod }k\} \vert.
$$
Let $\delta_t(k,l)$ be the number of solutions $x$ mod $k$ to the congruence
\begin{equation}
x^2g_t\equiv l \mbox{ mod } k. \label{23}
\end{equation}
Then it is easily seen that the condition (\ref{5}) with $S_t=S_t(Q)$ 
holds true for all positive $u\le Q/t$ and 
some absolute constant $C\ge 1$. Clearly, (\ref{6}) is also
valid. The remaining task is to bound $\delta_t(k,l)$.

If $(g_t,k)>1$, then $\delta_t(k,l)=0$ since $k$ and $l$ are supposed to be
coprime. Therefore, we can assume that $(g_t,k)=1$. Let $g$ mod $k$ 
be the multiplicative inverse of
$g_t$ mod $k$, i.e. $gg_t\equiv 1$ mod $k$. Put $l^*=gl$. 
Then $(\ref{23})$ is 
equivalent to
\begin{equation}
x^2\equiv l^* \mbox{ mod } k. \label{24}
\end{equation} 
Taking into account that $(k,l^*)=1$, and using some elementary facts on
the number 
of solutions of polynomial congruences modulo prime powers (see \cite{Qua}, 
for example), we see that  
(\ref{24}) has at most $2$ solutions if $k$ is a power of an odd prime and
at most $4$ solutions if $k$ is a power of 2. From this it follows that
for all $k\in \mathbbm{N}$ we have
$$
\delta_t(k,l)\le 2^{\omega(k)+1},
$$   
where $\omega(k)$ is the number of distinct 
prime divisors of $k$. For $k\le \sqrt{N}$
we have 
$$
2^{\omega(k)}\ll N^{\varepsilon}
$$
(see \cite{Han}). Therefore, (\ref{7}) holds with 
\begin{equation}
X:=c_6N^{\varepsilon}.\label{26}
\end{equation} 

Now we can apply Theorem 2. Combining (\ref{8}), (\ref{25}) and (\ref{26}), 
and taking into account that $Q\ge \sqrt{N}$, we obtain 
\begin{eqnarray}
\label{27}
& & \sum\limits_{q\in {\cal{S}}(Q)} 
\sum\limits_{\scriptsize \begin{array}{cccc} 
a=1\\(a,q)=1\end{array}}^q \left\vert S\left(\frac{a}{q}\right)\right\vert^2\\
&\le& c_7(\min\{QN^\varepsilon,N\}+Q)\left(\sqrt{N}\log\log 10 N+ \sqrt{Q}
\max\limits_{r\le\sqrt{N}} \sum\limits_{t\vert r} \frac{1}{f_t}\right)Z.
\nonumber
\end{eqnarray}

Next, we bound the function 
$$
\sigma(r)=\sum\limits_{t\vert r} \frac{1}{f_t}.
$$
Clearly, this function is multiplicative. If $r$ is a prime power $p^v$, then
$$
\sigma(r)\le 1+2\left(\frac{1}{p}+\frac{1}{p^2}+...\right)=1+\frac{2}{p-1}
\le \left(1+\frac{1}{p-1}\right)^2=\left(\frac{p^k}{\varphi\left(p^k\right)}
\right)^2.
$$
Hence, for all $r\in \mathbbm{N}$ we have
$$
\sigma(r)\le \left(\frac{r}{\varphi(r)}\right)^2\ll (\log \log 10r)^2.
$$
Therefore, from (\ref{27}), we obtain
\begin{eqnarray}
\label{28}
& & \sum\limits_{q\in {\cal{S}}(Q)} 
\sum\limits_{\scriptsize \begin{array}{cccc} 
a=1\\(a,q)=1\end{array}}^q \left\vert S\left(\frac{a}{q}\right)\right\vert^2\\
&\le& c_8(\log \log 10NQ)^2
(\min\{QN^\varepsilon,N\}+Q)\left(\sqrt{N}+ \sqrt{Q}\right)Z.
\nonumber
\end{eqnarray}
 
Considering the cases $QN^{\varepsilon}\le N$, $N<QN^{\varepsilon}\le 
N^{1+\varepsilon}$ and $N^{1+\varepsilon}<QN^{\varepsilon}$ separately,
we see that the right-hand side of (\ref{28}) is always bounded by
$$
\le c_9(\log\log 10NQ)^2\left(Q^{3/2}+N^{1/2+\varepsilon}Q\right)Z=:T(Q).
$$
Thus, we have 
\begin{eqnarray*}
\sum\limits_{q\in {\cal{S}}\cap(\sqrt{N},Q_1^2]} 
\sum\limits_{\scriptsize \begin{array}{cccc} 
a=1\\(a,q)=1\end{array}}^q
\left\vert S\left(\frac{a}{q}\right)\right\vert^2
&\le& \sum\limits_{j=0}^{\left[\log_2 (Q_1^2/\sqrt{N})\right]+1} 
T\left(\frac{Q_1^2}{2^j}\right)\\
&\ll& (\log\log 10NQ_1)^2\left(Q_1^3+N^{1/2+\varepsilon}Q_1^2\right)Z.
\end{eqnarray*}
Combining this with (\ref{22}), we obtain (\ref{9}). $\Box$\\
   
\section{Counting Farey fractions in short intervals}
Our later proof of Theorem 1 relies on the following variant of the large
sieve, which follows immediately from Theorem 2.11 in \cite{Lem}.\\ 

{\bf Lemma 1:} \begin{it} Let $\left(\alpha_r\right)_{r\in\mathbbm{N}}$ be a
sequence of real numbers and $\left(a_n\right)_{n\in\mathbbm{N}}$ be a
sequence of complex numbers. Define the trigonometrical polynomial
$S(\alpha)$ as in (\ref{0}). Suppose that $0<\Delta\le 1/2$ and 
$R\in \mathbbm{N}$. Put 
$$
K(\Delta):=\max\limits_{\alpha\in \mathbbm{R}} 
\sum\limits_{\scriptsize \begin{array}{cccc} r=1\\
\vert\vert \alpha_r -\alpha\vert\vert\le \Delta \end{array}}^R 1,
$$
where $\vert\vert x \vert\vert$ denotes the distance of a real $x$
to its closest integer.
Then 
$$
\sum\limits_{r=1}^R \left\vert S\left(\alpha_r\right)\right\vert^2
\le c_{10}K(\Delta)(N+\Delta^{-1})Z.
$$
\end{it}\\

In our situation, the sequence $\alpha_1,...,\alpha_R$ equals the sequence
of Farey fractions $a/q$ with $q\in {\cal{S}}$, $1\le a\le q$ and $(a,q)=1$. 
For $\alpha\in \mathbbm{R}$
put
$$
I(\alpha):=[\alpha-\Delta,\alpha+\Delta]
\ \ \mbox{ and }\ \ 
P(\alpha):= \sum\limits_{\scriptsize \begin{array}{cccc} 
q\in {\cal{S}}, (a,q)=1\\
a/q\in I(\alpha) \end{array}} 1.
$$
Then we have
\begin{equation}
K(\Delta)=\max\limits_{\alpha\in \mathbbm{R}} P(\alpha).\label{60}
\end{equation}
The next lemma provides an estimate for $P(\alpha)$.\\
       
{\bf Lemma 2:} \begin{it} Suppose that 
$0\le M\le Q$, ${\cal{S}}\subset (M,M+Q]$ and $0<\Delta\le 1/2$. 
Let $\alpha\in\mathbbm{R}$. 
Define $A_t(u,k,l)$ as in Theorem 1. Put
$$
U:=\left\{\begin{array}{llll} 1, & \mbox{ if }\ M< 1/\sqrt{\Delta},\\ \\
0, & \mbox{ otherwise.}\end{array}\right.
$$
Then,
\begin{eqnarray}
P(\alpha)&\le& c_{11}\left(U+\max\limits_{r\le 1/\sqrt{\Delta}}\ 
\max\limits_{\Delta\le z\le \sqrt{\Delta}/r} \max\limits_{\scriptsize
\begin{array}{cccc} h\in \mathbbm{Z}\\(h,r)=1\end{array}}\right. \label{10} \\
& & \left. \sum\limits_{t\vert r}
\sum\limits_{\scriptsize \begin{array}{cccc} 
0<m\le 4rzQ/t\\ (m,r/t)=1\end{array}} 
A_t\left(\frac{\Delta Q}{tz},\frac{r}{t},hm\right)\right).\nonumber
\end{eqnarray}
\end{it}\\

{\bf Proof:} We adapt Wolke's method used to prove Lemma 1 in \cite{Wol}. Let
\begin{equation}
\tau:=\frac{1}{\sqrt{\Delta}}.\label{P1}
\end{equation}
Then, by Dirichlet's approximation theorem, $\alpha$ can be written in the form
\begin{equation}
\alpha=\frac{b}{r}+z, \ \ \mbox{ where }\ \  r\le \tau,\ (b,r)=1,\ 
\vert z\vert \le \frac{1}{r\tau}.\label{P2}
\end{equation}
For $r\le \tau$ we have
\begin{equation}
\frac{1}{r\tau}\ge \frac{1}{\tau^2}=\Delta.\label{P3}
\end{equation}

We first note that we can restrict ourselves to the case when 
\begin{equation}
z\ge \Delta.\label{P4}
\end{equation}
If $\vert z\vert<\Delta$, then
$$ 
P(\alpha)\le P(b/r-\Delta)+P(b/r+\Delta).
$$ 
Furthermore, by (\ref{P3}), we have $\Delta\le 1/(r\tau)$. Therefore
this case can
be reduced to the case $\vert z\vert=\Delta$. 
Moreover, as $P(\alpha)=P(-\alpha)$, we can choose
$z$ positive. So we can assume (\ref{P4}).  

Write
\begin{equation}
\Delta_1:=\frac{\Delta}{z}\ \ \ \ \ \ (0<\Delta_1\le 1),\label{P5}
\end{equation}
\begin{equation}
j:=\left[\frac{1}{\Delta_1}\right]+1,\ \ \ \ \ \ 
\Delta_2:=\frac{1}{j}.\label{P6}
\end{equation}
Then, obviously,
\begin{equation}
\frac{\Delta_1}{2}\le \Delta_2\le \Delta_1.\label{P7}
\end{equation}
Put
\begin{equation}
y_i:=M+i\Delta_2Q,\ \ \ \ \ \  Y_i=(y_i,y_{i+1}],\ \ \ \ \ \ (i=0,1,...)
\label{P8}
\end{equation}
and 
\begin{equation}
P(\alpha)=\sum\limits_{i=0}^{j-1} \ \sum\limits_{q\in {\cal{S}}\cap Y_i}
\sum\limits_{\scriptsize \begin{array}{cccc} (a,q)=1\\ a/q\in I(\alpha)
\end{array}} 1 =: \sum\limits_{i=0}^{j-1} P_i(\alpha),\ \mbox{ say}.\label{P9}
\end{equation}

Now, for $0\le i\le j-1$ and 
\begin{equation}
q\in Y_i,\ \ \ \ \ \ (a,q)=1,\ \ \ \ \ \ \frac{a}{q}\in I(\alpha),\label{P10}
\end{equation}
we have $q(\alpha-\Delta)\le a\le q(\alpha+\Delta)$ or, by (\ref{P2}) and 
(\ref{P4}), $y_ir(z-\Delta)\le ar-bq\le y_{i+1}r(z+\Delta)$. Write
\begin{equation}
W_i:=[w_i,w_i^{\prime}]:=[y_ir(z-\Delta),y_{i+1}r(z+\Delta)].\label{Z3}
\end{equation}
Thus, (\ref{P10}) implies $ar-bq=m$, $m\in W_i$. If $m=0$, then 
$q=r$ since $(a,q)=1=(b,r)$. Hence,
\begin{equation}
P_i(\alpha)\le \nu_i+\sum\limits_{\scriptsize \begin{array}{cccc}
m\in W_i\\ m\not=0\end{array}}
\sum\limits_{\scriptsize 
\begin{array}{cccc} q\in {\cal{S}}\cap 
Y_i\\ q \equiv -\overline{b}m\mbox{ mod } r
\end{array}} 1, \label{P11}
\end{equation}
where 
$$
\nu_i:=\left\{\begin{array}{llll} 1, & \mbox{ if } \ r\in Y_i,\\ \\
0, & \mbox{ otherwise,}\end{array}\right.
$$
and $\overline{b}$ mod $r$ is 
the multiplicative inverse of $b$ mod $r$, {\it i.e.} $\overline{b}b\equiv
1$ mod $r$. 

For $t\in \mathbbm{N}$ write
$$
W_i/t:=[w_i/t,w_i^{\prime}/t],\ \ \ \ \ \  Y_i/t:=(y_i/t,y_{i+1}/t].
$$
Then from  (\ref{P11}) it follows that
\begin{equation}
P_i(\alpha)\le \nu_i+\sum\limits_{t\vert r} 
\sum\limits_{\scriptsize \begin{array}{cccc} 
m\in W_i/t\\ (m,r/t)=1\\ m\not=0\end{array}}\sum\limits_{\scriptsize 
\begin{array}{cccc} q\in {\cal{S}}_t\cap 
(Y_i/t)\\ q \equiv -\overline{b}m\mbox{ mod } r/t
\end{array}} 1.\label{P12}
\end{equation}
By (\ref{P5}), (\ref{P7}) and (\ref{P8}), we get
$y_{i+1}-y_i\le \Delta Q/z$.
Thus, from (\ref{P12}), we obtain
\begin{equation}
P_i(\alpha)\le \nu_i+\sum\limits_{t\vert r} 
\sum\limits_{\scriptsize \begin{array}{cccc} 
m\in W_i/t\\ (m,r/t)=1\\ m\not= 0\end{array}} 
A_t\left(\frac{\Delta Q}{tz},\frac{r}{t},-\overline{b}m\right).\label{P13}
\end{equation}

From (\ref{P5}), (\ref{P6}), (\ref{P7}), (\ref{P8}) and (\ref{Z3})
it follows that 
$w_i^{\prime}\le rz(y_{i+1}+2Q\Delta_2)=rzy_{i+3}$
and, similarly, $w_i\ge rzy_{i-2}$.
Moreover, we have $w_0\ge 0$ by (\ref{P4}) and $M\ge 0$. 
Thus, the intervals $W_i$ $(0\le i\le j-1)$ 
cover the interval $V=[0,rzy_{j+2}]$
at most 6 times. Furthermore, by (\ref{P6}), (\ref{P8}) and $M\le Q$, 
we have $rzy_{j+2}\le 4rzQ$. 
Hence, (\ref{P9}) and (\ref{P13}) imply 
\begin{equation}
P(\alpha)\le U+6 \sum\limits_{t\vert r} 
\sum\limits_{\scriptsize \begin{array}{cccc} 
0<m\le 4rzQ/t\\ (m,r/t)=1\end{array}} 
A_t\left(\frac{\Delta Q}{tz},\frac{r}{t},-\overline{b}m\right),\label{P14}
\end{equation}
where $U$ is defined as in Lemma 2.
From (\ref{P1}), (\ref{P2}), (\ref{P4}) and (\ref{P14}), we deduce 
(\ref{10}). This completes the proof. $\Box$\\

\section{Proofs of Theorems 1 and 2}
We are now in a position to prove our Theorems 1 and 2.\\ 
 
{\bf Proof of Theorem 1:} Taking $\Delta:=1/N$, the result of Theorem 1
follows from Lemma 1, (\ref{60}) and Lemma 2. $\Box$\\

{\bf Proof of Theorem 2:} 
Suppose that $r\le \sqrt{N}$, $1/N\le z\le 1/(r\sqrt{N})$, $(h,r)=1$ and
$t\vert r$.
Then,
by condition (\ref{5}), we have
\begin{equation}
\sum\limits_{\scriptsize \begin{array}{cccc} 
0<m\le 4rzQ/t\\ (m,r/t)=1\end{array}} 
A_t\left(\frac{Q}{tzN},\frac{r}{t},hm\right)
\le C\left(1+\frac{tS_t}{rzN}\right) F_h(r,t),\label{A1}
\end{equation}
where 
$$
F_h(r,t):=\sum\limits_{\scriptsize \begin{array}{cccc} 
0<m\le 4rzQ/t\\ (m,r/t)=1\end{array}}
\delta_t(r/t,hm).
$$

If $z\ge 1/(4Q)$, then the condition (\ref{6}) gives
\begin{equation}
F_h(r,t) \le \frac{8rzQ}{t}. \label{A3}
\end{equation}
If $z<1/(4Q)$, then the conditions (\ref{6}) and (\ref{7}) imply
\begin{equation}
F_h(r,t) \le \min\left\{\frac{4rzQX}{t},\frac{r}{t}\right\}. \label{A4}
\end{equation}  

For $r\le \sqrt{N}$ and $1/N\le z\le 1/(r\sqrt{N})$,   
from (\ref{A1}), (\ref{A3}) and (\ref{A4}), we derive
\begin{eqnarray}
& & N\ \sum\limits_{t\vert r}  
\sum\limits_{\scriptsize \begin{array}{cccc} 
0<m\le 4rzQ/t\\ (m,r/t)=1\end{array}} 
A_t\left(\frac{Q}{tzN},\frac{r}{t},hm\right) \label{A5}\\
&\le& c_{12}C\left(\min\left\{QX,N\right\}+Q\right)\left(
\sqrt{N}\sum\limits_{t\vert r}\frac{1}{t}+\sum\limits_{t\vert r} S_t\right)
\nonumber
\end{eqnarray} 
by a short calculation. From Theorem 1, (\ref{A5}) and 
$$
\sum\limits_{t\vert r}\frac{1}{t}\le \prod\limits_{p\vert r}
\left(1+\frac{1}{p}+\frac{1}{p^2}+...\right) = \prod\limits_{p\vert r} 
\frac{p}{p-1} =
\frac{r}{\varphi(r)}\le c_{13}\log\log 10r,
$$
we obtain the result of Theorem 2. $\Box$  \\ \\ \\
{\bf Acknowledgement.} This paper was written when the author held a
postdoctoral position at the Harish-Chandra Research Institute at Allahabad 
(India). The author wishes to thank this institute for financial support.\\

\end{document}